\documentclass[11pt]{article}

\usepackage[T1]{fontenc}
\usepackage[utf8]{inputenc}
\usepackage{amsmath,amssymb,epsfig,bbm,mathabx,textcomp}
\usepackage{stmaryrd}
\usepackage[colorlinks]{hyperref}
\usepackage{comment}
\usepackage{color}

\usepackage[acronym,toc]{glossaries}
\usepackage{glossary-longragged}
\usepackage{glossary-superragged}
\newglossary[slg]{symbols}{sym}{sbl}{List of Symbols}

\usepackage[textsize=small]{todonotes}
\usepackage{varwidth}

\usepackage{url}

\usepackage{enumitem}
\setlist[itemize]{labelindent=*, leftmargin=.5 truecm,nosep}

\usepackage[nottoc]{tocbibind}

\usepackage{fancyhdr}
\usepackage[us,12hr]{datetime} 
\fancypagestyle{plain}{
\fancyhf{}
\rfoot{ {\ddmmyyyydate\today}}
\lfoot{\thepage}
}

\usepackage{amsthm}

\theoremstyle{plain}
\newtheorem{defi}{Definition}[section]
\newtheorem{prop}[defi]{Proposition}
\newtheorem{theo}[defi]{Theorem}
\newtheorem{conj}[defi]{Conjecture}
\newtheorem{lemm}[defi]{Lemma}
\newtheorem{coro}[defi]{Corollary}

\theoremstyle{definition}
\newtheorem{rema}[defi]{Remark}
\newtheorem{exem}[defi]{Example}
\newtheorem{exems}[defi]{Examples}

\newcommand{\bdefi}{\begin{defi}}
\newcommand{\edefi}{\end{defi}}
\newcommand{\bprop}{\begin{prop}}
\newcommand{\eprop}{\end{prop}}
\newcommand{\btheo}{\begin{theo}}
\newcommand{\etheo}{\end{theo}}
\newcommand{\blemm}{\begin{lemm}}
\newcommand{\brema}{\begin{rema}}
\newcommand{\erema}{\end{rema}}
\newcommand{\bexer}{\begin{exem}}
\newcommand{\eexer}{\end{exem}}
\newcommand{\bexem}{\begin{exem}}
\newcommand{\eexem}{\end{exem}}
\newcommand{\bexems}{\begin{exems}}
\newcommand{\eexems}{\end{exems}}
\newcommand{\bconj}{\begin{conj}}
\newcommand{\econj}{\end{conj}}
\newcommand{\elemm}{\end{lemm}}
\newcommand{\bcoro}{\begin{coro}}
\newcommand{\ecoro}{\end{coro}}

\usepackage{mathrsfs}
\renewcommand\mathcal{\mathscr}

\newcommand{\A}{{\cal A}}

\newcommand{\C}{{\cal C}}

\newcommand{\G}{{\cal G}}

\newcommand{\M}{{\cal M}}

\newcommand{\OOO}{{\cal O}}
\renewcommand{\P}{{\cal P}}
\newcommand{\Q}{{\cal Q}}

\newcommand{\T}{{\cal T}}

\newcommand{\maths}[1]{{\mathbb #1}}  

\newcommand{\CC}{\maths{C}}

\newcommand{\HH}{\maths{H}}
\newcommand{\JJ}{\maths{J}}

\newcommand{\NN}{\maths{N}}

\newcommand{\PP}{\maths{P}}
\newcommand{\QQ}{\maths{Q}}
\newcommand{\RR}{\maths{R}}

\newcommand{\ZZ}{\maths{Z}}

\renewcommand{\lll}{{\mathfrak l}}

\newcommand{\sss}{{\mathfrak s}}

\newcommand{\uG}{\underline{G}}
\newcommand{\uH}{\underline{H}}

\newcommand{\bs}{\backslash}
\newcommand{\ga}{\gamma}
\newcommand{\Ga}{\Gamma}
\newcommand{\ov}[1]{{\overline{#1}}} 
\newcommand{\ra}{\rightarrow}

\newcommand{\transl}[1]{{\mathfrak t}_{#1}}
\newcommand{\dil}[1]{{\mathfrak h}_{#1}}

\newcommand{\Ad}{\operatorname{Ad}}

\newcommand{\cen}{\operatorname{cen}}

\newcommand{\Gal}{\operatorname{Gal}}
\newcommand{\gengeod}
{\operatorname{\widecheck{\G\,}\!\!}}

\newcommand{\hahi}
{\operatorname{AHI}}
\newcommand{\phahi}
{\operatorname
{\PP\!\hahi}}
\newcommand{\uphahi}{\underline{\phahi}}

\newcommand{\Heis}{\operatorname{Heis}}
\newcommand{\HS}{\mathcal{H\!S}}

\renewcommand{\Im}{{\operatorname{Im}}}

\newcommand{\mat}{\M}

\newcommand{\n}{\operatorname{\tt n}}

\newcommand{\ray}{\operatorname{rad}}
\renewcommand{\Re}{{\operatorname{Re}}}

\newcommand{\stab}{\operatorname{Stab}}
\newcommand{\tr}{\operatorname{\tt tr}}
\newcommand{\Tr}{\operatorname{Tr}}

\newcommand{\hdr}{{\HH}^2_\RR}

\newcommand{\hdc}{{\HH}^2_\CC}

\newcommand{\GL}{\operatorname{GL}}
\newcommand{\PGL}{\operatorname{PGL}}

\newcommand{\PSL}{\operatorname{PSL}}
\newcommand{\SO}{\operatorname{SO}}
\newcommand{\PO}{\operatorname{PO}}

\newcommand{\PU}{\operatorname{PU}}

\newcommand{\PSLC}{\operatorname{PSL}_{2}(\CC)}

\newcommand{\PSLR}{\operatorname{PSL}_{2}(\RR)}

\usepackage{mathtools}
\makeatletter
\DeclareRobustCommand\widecheck[1]{{\mathpalette\@widecheck{#1}}}
\def\@widecheck#1#2{%
    \setbox\z@\hbox{\m@th$#1#2$}%
    \setbox\tw@\hbox{\m@th$#1%
       \widehat{%
          \vrule\@width\z@\@height\ht\z@
          \vrule\@height\z@\@width\wd\z@}$}%
    \dp\tw@-\ht\z@
    \@tempdima\ht\z@ \advance\@tempdima2\ht\tw@ \divide\@tempdima\thr@@
    \setbox\tw@\hbox{%
       \raise\@tempdima\hbox{\scalebox{1}[-1]{\lower\@tempdima\box
\tw@}}}%
    {\ooalign{\box\tw@ \cr \box\z@}}}
\makeatother

\newcounter{fig}

\title{A classification of $\RR$-Fuchsian subgroups \\
of Picard modular groups}
\author{Jouni Parkkonen \and Fr\'ed\'eric Paulin} 

\begin{document}

\bibliographystyle{../alphanum}
\maketitle
\begin{abstract} 
    Given an imaginary quadratic extension $K$ of $\QQ$, we classify
    the maximal nonelementary subgroups of the Picard modular group
    $\PU(1,2;\OOO_K)$ preserving a totally real totally geodesic
    plane in the complex hyperbolic plane $\HH^2_\CC$. We prove that
    these maximal $\RR$-Fuchsian subgroups are arithmetic, and
    describe the quaternion algebras from which they arise. For
    instance, if the radius $\Delta$ of the corresponding $\RR$-circle
    lies in $\NN-\{0\}$, then the stabiliser arises from the
    quaternion algebra $\Big(\!\begin{array}{c} \Delta\,,\, |D_K|
    \\\hline\QQ\end{array} \!\Big)$.  We thus prove the existence of
    infinitely many orbits of $K$-arithmetic $\RR$-circles in the
    hypersphere of $\PP_2(\CC)$.
\end{abstract}


\section{Introduction}
\label{sec:intro}
Let $h$ be a Hermitian form with signature $(1,2)$ on $\CC^3$. The
projective unitary Lie group $\PU(1,2)$ of $h$ contains exactly two
conjugacy classes of connected Lie subgroups locally isomorphic to
$\PSLR$. The subgroups in one class are conjugate to
$\operatorname{P}(\operatorname{SU}(1,1)\times\{1\})$ and they
preserve a complex projective line for the projective action of
$\PU(1,2)$ on the projective plane $\PP_2(\CC)$, and those of the
other class are conjugate to $\PO(1,2)$ and preserve a maximal totally
real subspace of $\PP_2(\CC)$.  The groups $\PSLR$ and $\PU(1,2)$ act
as the groups of holomorphic isometries, respectively, on the upper
halfplane model $\hdr$ of the real hyperbolic space and on the
projective model $\hdc$ of the complex hyperbolic plane defined using
the form $h$.

If $\Ga$ is a discrete subgroup of $\PU(1,2)$, the intersections of
$\Ga$ with the connected Lie subgroups locally isomorphic to $\PSLR$
are its {\em Fuchsian subgroups}. The Fuchsian subgroups preserving a
complex projective line are called {\it $\CC$-Fuchsian}, and the ones
preserving a maximal totally real subspace are called {\it
  $\RR$-Fuchsian}. In \cite{ParPau17MS}, we gave a classification of
the maximal $\CC$-Fuchsian subgroups of the Picard modular groups, and
we explicited their arithmetic structures, completing work of
Chinburg-Stover (see Theorem 2.2 in version 3 of \cite{ChiSto11} and
\cite[Theo.~4.1]{ChiSto18}) and Möller-Toledo in \cite{MolTol15}, in
analogy with the result of Maclachlan-Reid \cite[Thm.~9.6.3]{MacRei03}
for the Bianchi subgroups in $\PSLC$. In this paper, we prove
analogous results for $\RR$-Fuchsian subgroups, thus completing an
arithmetic description of all Fuchsian subgroups of the Picard modular
groups. The classification here is more involved, as in some sense,
there are more $\RR$-Fuchsian subgroups than $\CC$-Fuchsian ones. Our
approach is elementary, some of the results can surely be obtained by
more sophisticated tools from the theory of algebraic groups.

Let $K$ be an imaginary quadratic number field, with discriminant
$D_K$ and ring of integers $\OOO_K$. We consider the Hermitian form
$h$ defined by
$$
(z_0,z_1,z_2)\mapsto -\frac{1}{2}\,z_0\,\overline{z_2}
-\frac{1}{2}\,z_2\,\overline{z_0}+ z_1\overline{z_1}\;.
$$ The {\it Picard modular group} $\Ga_K=\PU(1,2)\cap\PGL_3(\OOO_K)$
is a nonuniform arithmetic lattice of $\PU(1,2)$.\footnote{See for
  instance \cite[Chap.~5]{Holzapfel98} and subsequent works of Falbel,
  Parker, Francsics, Lax, Xie, Wang, Jiang, Zhao and many others, for
  information on these groups, using different Hermitian forms of
  signature (2,1) defined over $K$.}  In this paper, we classify the
maximal $\RR$-Fuchsian subgroups of $\Ga_K$, and we explicit their
arithmetic structures. The results stated in this introduction do not
depend on the choice of the Hermitian form $h$ of signature $(2,1)$
defined over $K$, since the algebraic groups over $\QQ$ whose groups
of $\QQ$-points are $\PU(1,2)\cap \PGL_3(K)$ depend up to
$\QQ$-isomorphism only on $K$ and not on $h$, see for instance
\cite[\S~3.1]{Stover11}, so that the Picard modular group $\Ga_K$ is
well defined up to commensurability.


Let $I_3$ be the identity matrix and let $I_{1,2}$ be the matrix of $h$. Let
$$
\hahi(\QQ)=\{Y\in\M_3(K)\;:\;\; 
Y^*I_{1,2}Y=I_{1,2} \;\;{\rm and}\;\;Y\,\ov{Y}=I_3\}
$$ 
be the set of $\QQ$-points of an algebraic subset defined over
$\QQ$, whose real points consist of the matrices of the Hermitian
anti-holomorphic linear involutions $z\mapsto Y\,\ov{z}$ of
$\CC^3$. For instance,
$$
Y_\Delta=\begin{pmatrix} 0 & 0 & \frac{1}{\overline{\Delta}}\\ 0 & 1 &
0 \\ \Delta & 0 & 0
\end{pmatrix}
$$  
belongs to $\hahi(\QQ)$ for every $\Delta\in K^\times$. The group
$\operatorname{U}(1,2)$ acts transitively on $\hahi(\RR)$ by
$$
(X,Y)\mapsto X\,Y\,\ov{X}^{\,-1}
$$ 
for all $X\in \operatorname{U}(1,2)$ and $Y\in\hahi(\RR)$.  In Section
\ref{sect:class}, we prove the following result that describes the
collection of maximal $\RR$-Fuchsian subgroups of the Picard modular
groups $\Ga_K$.

\btheo\label{theo:classification} The stabilisers in $\Ga_K$ of the
projectivized rational points in $\hahi(\QQ)$ are arithmetic maximal
nonelementary $\RR$-Fuchsian subgroups of $\Ga_K$.  Every maximal
nonelementary $\RR$-Fuchsian subgroup of $\Ga_K$ is commensurable up
to conjugacy in $\PU(1,2)\cap\PGL_3(K)$ with the stabiliser
$\Ga_{K,\,\Delta}$ in $\Ga_K$ of the projective class of $Y_\Delta$,
for some $\Delta\in\OOO_K-\{0\}$.  
\etheo

A nonelementary $\RR$-Fuchsian subgroup $\Ga$ of $\PU(1,2)$ {\it
  arises from} a quaternion algebra $\Q$ over $\QQ$ if $\Q$ splits
over $\RR$ and if there exists a Lie group epimorphism $\varphi$ from
$\Q(\RR)^1$ to the conjugate of $\PO(1,2)$ containing $\Ga$ such that
$\Ga$ and $\varphi(\Q(\ZZ)^1)$ are commensurable.  In Section
\ref{sect:ternquad}, we use the connection between quaternion algebras
and ternary quadratic forms to describe the quaternion algebras from
which the maximal nonelementary $\RR$-Fuchsian subgroups of the Picard
modular groups $\Ga_K$ arise.

\btheo \label{theo:main} 
For every $\Delta\in \OOO_K-\{0\}$, the maximal nonelementary
$\RR$-Fuchsian subgroup $\Ga_{K,\,\Delta}$ of $\Ga_K$ arises from the
quaternion algebra with Hilbert symbol 
$\big(\frac{2 \Tr_{K/\QQ}\Delta
  ,\, N_{K/\QQ}(\Delta) \,|D_K|}{\QQ}\big)$ if $\Tr_{K/\QQ}\Delta\ne 0$
and from $\big(\frac{1,\,1}{\QQ}\big)\simeq \M_2(\QQ)$ otherwise.
\etheo

This arithmetic description has the following geometric consequence.
Recall that an {\it $\RR$-circle} is a topological circle which is the
intersection of the {\it Poincaré hypersphere}
$$
\HS=\{[z]\in\PP_2(\CC)\;:\;h(z)=0\}
$$ 
with a maximal totally real subspace of $\PP_2(\CC)$. It is {\it
  $K$-arithmetic} if its stabiliser in $\Ga_K$ has a dense orbit in
it.

\bcoro \label{coro:intro} There are infinitely many $\Ga_K$-orbits of
$K$-arithmetic $\RR$-circles in the hypersphere $\HS$.  
\ecoro

The figure below shows the image under vertical projection from
$\partial_\infty\HH^2_\CC$ to $\CC$ of part of the
$\Ga_{\QQ(i)}$-orbit of the standard infinite $\RR$-circle, which is
$\QQ(i)$-arithmetic. The image of each finite $\RR$-circle is a
lemniscate. We refer to Section \ref{sec:R-circles} and \cite[\S
  4.4]{Goldman99} for an explanation of the terminology. See the main
body of the text for other pictures of $K$-arithmetic $\RR$-circles.

\medskip 

\begin{center}
\includegraphics[width=9cm]{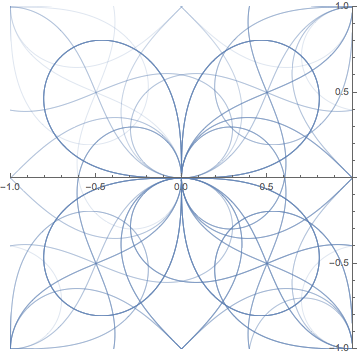}
\end{center}

\section{The complex hyperbolic plane}
\label{sec:cxhyp}


Let $h$ be the nondegenerate Hermitian form on $\CC^3$ defined by
$$
h(z)=z^*I_{1,2}z=-\Re(z_0\overline{z_2})+|z_1|^2\,,
$$
where $I_{1,2}$ is the antidiagonal matrix
$$
I_{1,2}=\begin{pmatrix} \ \ 0 & 0 & -\frac{1}{2} \\ \ \ 0 & 1 &\ \  0 \\ 
-\frac{1}{2} & 0 & \ \ 0 \end{pmatrix}\,.
$$ 
A point $z=(z_0,z_1,z_2)\in\CC^3$ and the corresponding element
$[z]=[z_0:z_1:z_2]\in \PP_2(\CC)$ (using homogeneous coordinates) is
{\em negative, null or positive} according to whether $h(z)<0$,
$h(z)=0$ or $h(z)>0$.  The {\it negative/null/positive cone} of $h$ is
the subset of negative/null/positive elements of $ \PP_2(\CC)$.

The negative cone of $h$ endowed with the distance $d$ defined by
$$
\cosh^2d([z],[w])=\frac{|\langle z,w\rangle|^2}{h(z)\,h(w)}\;,
$$
where $\langle \cdot\,,\cdot\rangle$ is the sesquililnear form
associated with $h$, is the {\em complex hyperbolic plane} $\hdc$. The
distance $d$ is the distance of a Riemannian metric with pinched
negative sectional curvature $-4\le K\le -1$.  The null cone of $h$ is
the Poincar\'e hypersphere $\HS$, which is naturally identified
with the boundary at infinity of $\hdc$.

The Hermitian form $h$ in this paper differs slightly from the one we
used in \cite{ParPau10GT,ParPau17MA,ParPau17MS} and from the main
Hermitian form used by Goldman and Parker (see \cite{Goldman99,
  Parker12, Parker16}).  Hence we will need to give some elementary
computations that cannot be found in the literature. This form is a
bit more appropriate for arithmetic purposes concerning $\RR$-Fuchsian
subgroups, as it allows us to consider $\ZZ$-points of our linear
algebraic groups and not their $2\ZZ$-points.

Let $\operatorname{U}(1,2)$ be the linear group of $3\times3$
invertible matrices with complex coefficients preserving the Hermitian
form $h$. Let $\PU(1,2)=\operatorname{U}(1,2)/ \operatorname{U}(1)$ be
its associated projective group, where $\operatorname{U}(1)=
\{\zeta\in\CC \;:\;|\zeta|=1\}$ acts by scalar multiplication. We
denote by $[X]=[a_{ij}]_{1\leq i,j\leq n}\in\PU(1,2)$ the image of
$X=(a_{ij})_{1\leq i,j\leq n}\in\operatorname{U}(1,2)$. The linear
action of $\operatorname{U}(1,2)$ on $\CC^3$ induces a projective
action of $\PU(1,2)$ on $\PP_2(\CC)$ that preserves the negative, null
and positive cones of $h$ in $\PP_2(\CC)$, and is transitive on each
of them.

If
$$
X=\begin{pmatrix} a & \ov{\ga} & b\\
    \alpha & A & \beta \\ c & \ov{\delta} & d\end{pmatrix}\in\M_3(\CC),
\textrm{ then }\quad
I_{1,2}^{-1}X^*I_{1,2}=\begin{pmatrix} \ \ \ov{d} & -2\ov\beta & 
\ \ \ov{b} \vspace{.1cm}\\ 
-\frac{\delta}{2} & \ \ \ov{A} & -\frac{\ga}{2} \vspace{.1cm}\\ 
\ \ \ov{c} & -2\ov{\alpha} & \ \ \ov{a}\end{pmatrix}.
$$
The matrix $X$ belongs to $\operatorname{U}(1,2)$ if and only if $X$ is
invertible with inverse $I_{1,2}^{-1}X^*I_{1,2}$, that is, if and only if
\begin{equation}\label{eq:equationsU}
\left\{\begin{array}{l}
a\ov{d}+b\overline{c}-\frac{1}{2}\delta\ov{\ga} =1\\
\ov{d}\alpha+\overline{c}\beta-\frac{1}{2}A\delta =0\\
c\ov{d}+d\ov{c}-\frac{1}{2}|\delta|^2 = 0 \\ 
A \ov{A}-2\alpha\ov{\beta}- 2\beta\ov{\alpha} =1 \\
a\ov{b}+b\overline{a}-\frac{1}{2}|\ga|^2 = 0 \\ 
\ov{b}\alpha+\ov{a}\beta-\frac{1}{2}A\ga =0 \;.
\end{array}\right.
\end{equation}

\brema\label{rem:cnultriangsup} A matrix $X\in \operatorname{U}(1,2)$
in the above form is upper triangular if and only if $c=0$. Indeed,
then the third equality in Equation \eqref{eq:equationsU} implies that
$\delta=0$. The first two equations then become $a\ov{d}=1$ and
$\ov{d}\alpha=0$, so that $\alpha =0$.  \erema

The
{\em Heisenberg group}
$$
\Heis_3=\big\{[w_0:w:1]\in\PP_2(\CC):\Re\, w_0=|w|^2\big\} 
$$ 
with law $[w_0:w:1][w'_0,w':1]= [w_0+w'_0+2w'\,\overline{w},
  w+w':1]$ is identified with $\CC\times\RR$ by the coordinate mapping
$[w_0:w:1]\mapsto(w,\Im\, w_0)=(\zeta,v)$.  It acts isometrically on
$\hdc$ and simply transitively on $\HS- \{[1:0:0]\}$ by {\em
  Heisenberg translations}
$$
\transl{\zeta,v}=
\begin{bmatrix} 1 & 2\,\ov{\zeta} & |\zeta|^2+iv\\ 0 & 1 &
\zeta \\ 0 & 0 & 1\end{bmatrix}\in\PU(2,1)
$$ 
with $\zeta\in\CC$ and $v\in\RR$. Note that $\transl{\zeta,v}^{-1}
=\transl{-\zeta,-v}$ and $\overline{\transl{\zeta,v}}=
\transl{\overline{\zeta},-v}$.  The {\it Heisenberg dilation} with
factor $\lambda\in\CC^\times$ is the element 
$$
\dil\lambda=\begin{bmatrix} \lambda & 0 & 0\\ 0 & 1 & 0 \\ 0
& 0 & \frac{1}{\ov{\lambda}}\end{bmatrix}\in\PU(1,2)\,,
$$
which normalizes the group of Heisenberg translations.
The subgroup of $\PU(1,2)$ generated by Heisenberg translations and
Heisenberg dilations is called the group of {\it Heisenberg
  similarities}.

We end this subsection by defining the discrete subgroup of $\PU(1,2)$
whose $\RR$-Fuchsian subgroups we study in this paper.

Let $K$ be an imaginary quadratic number field, with $D_K$ its
discriminant, $\OOO_K$ its ring of integers, $\Tr:z\mapsto
z+\overline{z}$ its trace and $N:z\mapsto|z|^2=z\,\overline{z}$ its
norm.  Recall\footnote{See for instance \cite{Samuel67}.} that there
exists a squarefree positive integer $d$ such that $K=\QQ(i\sqrt{d})$,
that $D_K=-d$ and $\OOO_K=\ZZ[\frac{1+i\sqrt{d}}{2}]$ if $d\equiv
-1\mod 4$, and that $D_K=-4d$ and $\OOO_K=\ZZ[i\sqrt{d}]$
otherwise. Note that $\OOO_K$ is stable under conjugation, and that
$\Tr$ and $N$ take integral values on $\OOO_K$. A {\it unit} in
$\OOO_K$ is an invertible element in $\OOO_K$. Since
$N:K^\times\ra\RR^\times$ is a group morphism, we have $N(x)=1$ for
every unit $x$ in $\OOO_K$.

The {\em Picard modular group}
$$
\Ga_K=\PU(1,2;\OOO_K)=\PU(1,2)\cap \PGL_3(\OOO_K)
$$
is a nonuniform lattice in $\PU(1,2)$.

\section{The space of $\RR$-circles}
\label{sec:R-circles}

A (maximal) {\em totally real subspace} $V$ of the Hermitian vector
space $(\CC^3,h)$ is the fixed point set of a Hermitian
antiholomorphic linear involution of $\CC^3$, or, equivalently, a
$3$-dimensional real linear subspace of $\CC^3$ such that $V$ and $\JJ
V$ are orthogonal, where $\JJ:\CC^3\to\CC^3$ is the componentwise
multiplication by $i$. The intersection with $\HH^2_\CC$ of the image
under projectivization in $\PP_2(\CC)$ of a totally real subspace is
called an $\RR${\em -plane} in $\HH^2_\CC$. The group $\PU(1,2)$ acts
transitively on the set of $\RR$-planes, the stabiliser of each
$\RR$-plane being a conjugate of $\PO(1,2)$.  Note that $\PO(1,2)$ is
equal to its normaliser in $\PU(1,2)$.

An {\it $\RR$-circle} $C$ is the boundary at infinity of an
$\RR$-plane. See \cite{Mostow73}, \cite[\S 4.4]{Goldman99} and
\cite[\S 9]{Jacobowitz90} for references on $\RR$-circles (introduced
by E.~Cartan).  An $\RR$-circle is {\it infinite} if it contains
$\infty=[1:0:0]$ and {\it finite} otherwise. The group of Heisenberg
similarities acts transitively on the set of finite $\RR$-circles and
on the set of infinite $\RR$-circles.

The {\it standard infinite $\RR$-circle} is
$$
C_\infty=\big\{[x_0:x_1:x_2]\;:\; x_0,x_1,x_2\in\RR,\;x_1^2-x_0x_2=0\big\}\;,
$$ 
which is the boundary at infinity of the intersection with
$\HH^2_\CC$ of the image in $\PP_2(\CC)$ of $\RR^3\subset \CC^3$. For
every $D\in\CC^\times$, the set
$$
C_{D}=\big\{[z_0:x_1:D\,\ov{z_0}]\;:\; z_0\in\CC,\;x_1\in\RR,\;
\;x_1^2-\Re(\,\ov{D}z_0^2)=0\big\}
$$ 
is a finite $\RR$-circle, which is the boundary at infinity of the
intersection with $\HH^2_\CC$ of the fixed point set of the projective
Hermitian anti-holomorphic involution
$$
[z_0:z_1:z_2]\mapsto [\frac{\ov{z_2}}{\ov{D}}:\ov{z_1}:D\,\ov{z_0}]\,.
$$ 
We call $C_1$ the {\it standard finite $\RR$-circle}.

Let $C$ be a finite $\RR$-circle. The {\it center} $\cen(C)$ of $C$ is
the image of $\infty=[1:0:0]$ by the unique projective Hermitian
anti-holomorphic involution fixing $C$. The {\it radius} $\ray(C)$ of
$C$ is $\lambda^2$ where $\lambda\in\CC^\times$ is such that there
exists a Heisenberg translation $\transl{}$ mapping $0=[0:0:1]$ to the
center of $C$ with $C=\transl{}\circ \dil{\lambda}(C_1)$.  For
instance, $\cen(C_D)=0$ and $\ray(C_D)= \frac{1}{\ov{D}}$, since the
Heisenberg dilations preserve $0$ and $C_D=
\dil{\frac{1}{\sqrt{\ov{D}}}} (C_1)$.  For every Heisenberg
translation $\transl{}$, we have $\cen(\transl{} C)=\transl{}\cen(C)$
and $\ray(\transl{} C)= \ray(C)$.  For every Heisenberg dilation
$\dil\lambda$, we have $\cen(\dil\lambda C)=\dil\lambda\cen(C)$ and
$\ray(\dil\lambda C)= \lambda^2\ray(C)$.

The image of a finite $\RR$-circle under the {\em vertical projection}
$(\zeta,v)\mapsto \zeta$ from $\Heis_3=\partial_\infty\HH^2\CC-
\{\infty\}$ to $\CC$ is a lemniscate, see \cite[\S 4.4.5]{Goldman99}.
The figure below shows on the left six images of the standard infinite
$\RR$-circle under transformations in $\Ga_{\QQ(\omega)}$ where
$\omega=\frac{-1+i\sqrt 3}2$ is the usual third root of unity, and on
the right their images in $\CC$ under the vertical projection.
\begin{center}
\includegraphics[width=6cm]{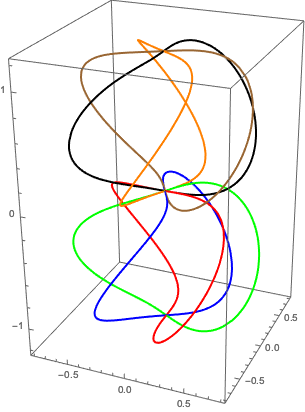}
\includegraphics[width=6cm]{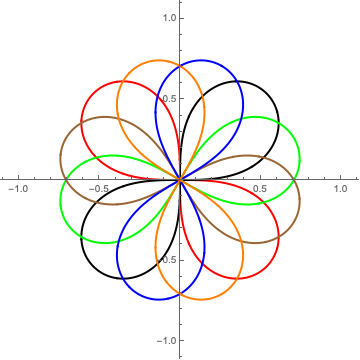}
\end{center}

Let us introduce more notation in order to describe the space of
$\RR$-circles, see \cite[\S 2.2.4]{Goldman99} for more background. A
$3\times 3$ matrix $Y$ with complex coefficients is called {\it
  unitary-symmetric} if it is Hermitian with respect to the Hermitian
form $h$ and invertible with inverse equal to its complex conjugate,
that is, if $Y^*I_{1,2}Y=I_{1,2}$ and $Y\,\ov{Y}=I_3$, where $I_3$ is
the $3\times 3$ identity matrix.  Note that for instance $I_3$ and,
for every $D\in\CC^\times$, the matrix
$$
Y_D=\begin{pmatrix} 
0 & 0 & \frac{1}{\overline{D}}\\ 0 & 1 & 0 \\ D & 0  & 0
\end{pmatrix}
$$ 
is unitary-symmetric. 

Let
$$
\hahi=\{Y\in\M_3(\CC)\;:\;\; 
Y^*I_{1,2}Y=I_{1,2} \;\;{\rm and}\;\;Y\,\ov{Y}=I_3\}
$$ 
be the set of unitary-symmetric matrices, which is a closed subset of
$\operatorname{U}(1,2)$, identified with the set of Hermitian
anti-holomorphic linear involutions $z\mapsto Y\,\ov{z}$ of
$\CC^3$. Note that $|\det Y|=1$ for any $Y\in\hahi$. Let
$$
\phahi=\{[Y]\in\PU(1,2)\;:\;\; Y\,\ov{Y}=I_3\}
$$ 
be the image of $\hahi$ in $\PU(1,2)$, that is, the quotient
$\operatorname{U}(1)\bs\hahi$ of $\hahi$ modulo scalar multiplications
by elements of $\operatorname{U}(1)$.  The group
$\operatorname{U}(1,2)$ acts transitively
on $\hahi$ by
$$
(X,Y)\mapsto X\,Y\,\ov{X}^{\,-1}
$$
for all $X\in \operatorname{U}(1,2)$ and $Y\in\hahi$, and the
stabiliser of $I_3$ is equal to $\operatorname{O}(1,2)$.

For every $Y\in\hahi$, we denote by $P_Y$ the intersection with
$\HH^2_\CC$ of the image in $\PP_2(\CC)$ of the set of fixed points of
$z\mapsto Y\ov{z}$. Note that $P_Y$ is an $\RR$-plane, which depends
only on the class $[Y]$ of $Y$ in $\PU(1,2)$. We denote by
$C_Y=\partial_\infty P_Y$ the $\RR$-circle at infinity of $P_Y$, which
depends only on $[Y]$. For instance, $C_\infty=C_{I_3}$ and
$C_D=C_{Y_D}$.

\newcommand\rcircles{\C_\RR}
\newcommand\rplanes{\P_\RR}

\medskip
Let $\rcircles$ be the set of $\RR$-circles, endowed with the topology
induced by the Hausdorff distance between compact subsets of
$\partial_\infty\HH^2_\CC$,\footnote{for any Riemannian distance on
  the smooth manifold $\partial_\infty\HH^2_\CC$} and let $\rplanes$
be the set of $\RR$-planes\footnote{which are closed subsets of
  $\HH^2_\CC$} endowed with the topology of the Hausdorff convergence
on compact subsets of $\HH^2_\CC$.

The projective action of $\PU(1,2)$ on the set of subsets of
$\PP_2(\CC)$ induces continuous transitive actions on $\rcircles$ and
$\rplanes$, with stabilisers of $C_\infty=C_{I_3}$ and $P_{I_3}$ equal
to $\PO(1,2)$.  We hence have a sequence of $\PU(1,2)$-equivariant
homeomorphisms
\begin{equation}\label{eq:homeospaceRcercle}
\begin{array}{ccccccc}
\vspace{.1cm} 
\PU(1,2)/\PO(1,2) & \longrightarrow & \phahi & \longrightarrow & 
\rplanes & \longrightarrow & \rcircles\\
\vspace{.1cm} 
{}[X]\PO(1,2) &\longmapsto & \big[X\ov{X}^{-1}\big]& & P &\longmapsto &
\partial_\infty P \;.\\ & & [Y]&\longmapsto & \ P_Y & &
\end{array}
\end{equation}

\blemm \label{lem:calccentrerayon}
Let $Y=\begin{pmatrix} a & \ov\ga & b\\ \alpha & A & \beta \\ 
c & \ov{\delta} & d \end{pmatrix} \in\hahi$.
\begin{enumerate}
\item[(1)] For every $[X]\in\PU(1,2)$, we have
$[X]\,C_Y= C_{XY\,\ov{X}^{\,-1}}$.
\item[(2)]  The $\RR$-circle $C_Y$ is infinite if and only if $c=0$.
\item[(3)] If the $\RR$-circle $C_Y$ is finite, then its center is 
$$
\cen(C_Y)=[Y]\,\infty=[a:\alpha:c]\;,
$$
and its radius is
$$
\ray(C_Y)=  \frac{\ov{A\,c}-\ov{\alpha}\,\delta}{\ov{c}^{\,2}} =
-\frac{c}{\ov{c}^{\,2}} \;\overline{\det Y}
\;.
$$
In particular, $\big|\ray(C_Y)\big|= |c|^{-1}$. 
\end{enumerate}
\elemm

\begin{proof} (1) This follows from the equivariance of the homeomorphisms in
Equation \eqref{eq:homeospaceRcercle}.

\medskip\noindent 
(2) Recall that $C_Y$ is the intersection with $\partial_\infty
\HH^2_\CC$ of the image in the projective plane of the set of fixed
points of the Hermitian anti-holomorphic linear involution $z\mapsto
Y\,\ov{z}$. Hence $\infty=[1:0:0]$ belongs to $C_Y$ if and only if the
image of $(1,0,0)$ by $Y$ is a multiple of $(1,0,0)$, that is, if and
only if $\alpha=c=0$. Using Remark \ref{rem:cnultriangsup}, this
proves the result.

\medskip\noindent 
(3) The first claim follows from the fact that the center of the
$\RR$-circle $C_Y$ is the image of $\infty=[1:0:0]$ under the
projective map associated with $z\mapsto Y\,\ov{z}$. In order to prove
the second claim, we start by the following lemma.

\blemm \label{lem:formnormhahi} For every $[Y]\in\phahi$, the center
of $C_Y$ is equal to $0=[0:0:1]$ if and only if there exists
$D\in\CC^\times$ such that $[Y]=[Y_D]$.  
\elemm

\begin{proof} 
We have already seen that $\cen(C_{Y_D})=\cen(C_{D})=0$. By the first
claim of Lemma \ref{lem:calccentrerayon} (3), if $\cen(C_Y)=0$, we have
$a=\alpha=0$. By the penultimate equality in Equation
\eqref{eq:equationsU}, we have $\ga=0$.  Since $Y\, \ov Y=I_3$, we
have $b\,\ov{c}=1$, $b\,\delta=0$, $b\,\ov{d}=0$ and
$\beta\,\ov{c}=0$, so that $Y=\begin{pmatrix} 0 & 0 &
\frac{1}{\ov{c}}\\ 0 & A & 0 \\ c & 0 & 0 \end{pmatrix}$ with
$|A|=1$. Since $[Y]=[\frac{1}{A} Y]$, the result follows with
$D=\frac{c}{A}$.  
\end{proof}

Now, let $\zeta=\frac{\alpha}{c}$, $v=\Im\,\frac{a}{c}$ and
$X=\begin{pmatrix} 1 & 2\,\overline{\zeta} & |\zeta|^2+iv\\ 0 & 1 &
\zeta \\ 0 & 0 & 1 \end{pmatrix}$. Note that since $Y\in
\operatorname{U}(1,2)$, we have
$$
|\alpha|^2-\Re(a\,\ov{c})=h(a,\alpha,c)=
h(Y(1,0,0)) = h(1,0,0)=0\,.
$$ 
Hence
$$
\Re\big(\frac{a}{c}\big)=\frac{1}{|c|^2}\,\Re(a\,\ov{c})=
\Big|\frac{\alpha}{c}\Big|^2=|\zeta|^2\;.
$$ 
The Heisenberg translation $\transl{\zeta,v}=[X]$ maps $0=[0:0:1]$
to $[\frac{a}{c}: \frac{\alpha}{c}:1] =\cen(C_Y)$. Since
$$
\cen(C_{X^{-1}Y\,\ov{X}}) =
\cen(\transl{\zeta,v}^{-1}C_Y) = \transl{\zeta,v}^{-1} \cen(C_Y) =0\,,
$$ 
and by Lemma \ref{lem:formnormhahi}, the element $X^{-1}Y\,\ov{X}
\in\hahi$ is anti-diagonal. A simple computation gives
$$
X^{-1}Y\,\ov{X}
= \begin{pmatrix} 0 & 0 & \frac{1}{\ov{c}}\\ 0 & A-\zeta \ov{\delta} &
  0 \\ c & 0 & 0 \end{pmatrix}.
$$ 
If $D= \frac{c}{A-\zeta \ov{\delta}}$, we hence have $[X^{-1}Y\ov{X}]
= [Y_D]$. Therefore
$$
\ray(C_Y)=\ray(\transl{\zeta,v}^{-1}C_Y)=\ray(C_{X^{-1}Y\ov{X}})=\ray(C_{Y_D})=
\frac{1}{\ov{D}}\;.
$$
Since $\det X=1$, we have $\det Y=-\,\frac{c}{\ov{c}}\,(A-\zeta
\ov{\delta})$, so that $D=-\,\frac{c^2}{\ov{c}\;\det Y}$.  The result
follows. 
\end{proof}

We end this section by describing the algebraic properties of the objects
in Equation \eqref{eq:homeospaceRcercle}.  We refer for instance to
\cite[\S 3.1]{Zimmer84} for an elementary introduction to algebraic
groups and their Zariski topology.

Let $\uG$ be the linear algebraic group defined over $\QQ$, with set
of $\RR$-points $\PU(1,2)$ and set of $\QQ$-points 
$$
\PU(1,2;K)=\PU(1,2)\cap\PGL_3(K)\;.
$$  
We identify $\uG$ with its image under the adjoint representation for
integral point purposes, so that $\uG(\ZZ)=\Ga_K$.

Since $I_{1,2}$ has rational coefficients, the set $\phahi$ of
unitary-symmetric matrices modulo scalars is the set of real points
$\phahi= \uphahi(\RR)$ of an affine algebraic subset $\uphahi$ defined
over $\QQ$ of $\uG$, whose set of rational points is 
$$
\uphahi(\QQ)=
\phahi\cap\,\uG(\QQ)= \phahi\cap\PGL_3(K)\;.
$$
The action of $\uG$ on $\uphahi$ defined by $([X],[Y])\mapsto
[X\,Y\,\ov{X}^{\,-1}]$ is algebraic defined over $\QQ$.  This notion
of rational point in $\phahi$ will be a key tool in the next section
in order to describe the maximal nonelementary $\RR$-Fuchsian
subgroups of $\Ga_K$.

\section{A description of the $\RR$-Fuchsian subgroups of $\Ga_K$}
\label{sect:class}

Our first result relates the nonelementary $\RR$-Fuchsian subgroups of
the Picard modular group $\Ga_K$ to the rational points in
$\phahi$. The proof of this statement is similar to the one of its
analog for $\CC$-Fuchsian subgroups in \cite{ParPau17MS}.

\bprop\label{prop:Fuchsianrationalpoint} The stabilisers in $\Ga_K$ of
the rational points in $\phahi$ are maximal nonelementary
$\RR$-Fuchsian subgroups of $\Ga_K$. Conversely, any maximal
nonelementary $\RR$-Fuchsian subgroup $\Ga$ of $\Ga_K$ fixes a unique
rational point in $\phahi$ and $\Ga$ is an arithmetic lattice in the
conjugate of $\PO(1,2)$ containing it.  
\eprop

\begin{proof}
Let $[Y]\in\uphahi(\QQ)$ be a rational point in $\phahi$.  Since
the action of $\uG$ on $\uphahi$ is algebraic defined over $\QQ$, the
stabiliser $\uH$ of $[Y]$ in $\uG$ is algebraic defined over $\QQ$.
Note that $\uH$ is semi-simple with set of real points a conjugate of
(the normaliser of $\PO(1,2)$ in $\PU(1,2)$, hence of)
$\PO(1,2)$. Therefore by the Borel-Harish-Chandra theorem
\cite[Thm.~7.8]{BorelHarishChandra62}, the group $\stab_{\Ga_K}[Y]=
\uH(\ZZ)$ is an arithmetic lattice in $\uH(\RR)$, and in particular is
a maximal nonelementary $\RR$-Fuchsian subgroup of $\Ga_K$.

Conversely, let $\Ga$ be a maximal nonelementary $\RR$-Fuchsian
subgroup of $\Ga_K$. Since it is nonelementary, its limit set
$\Lambda\Ga$ contains at least three points. Two $\RR$-circles having
three points in common are equal. Hence $\Ga$ preserves a unique
$\RR$-plane $P$. Let $Y\in\hahi$ be such that $P=P_Y$. By the
equivariance of the homeomorphisms in Equation
\eqref{eq:homeospaceRcercle}, $[Y]$ is the unique point in $\phahi$
fixed by $\Ga$.

Let $\uH$ be the stabiliser in $\uG$ of $[Y]$, which is a connected
algebraic subgroup of $\uG$ defined over $\RR$, whose set of real
points is conjugated to $\PO(1,2)$.  Since a nonelementary subgroup of
a connected algebraic group whose set of real points is isomorphic to
$\PSL_2(\RR)$ is Zariski-dense in it, and since the Zariski-closure of
a subgroup of $\uG(\ZZ)$ is defined over $\QQ$ (see for instance
\cite[Prop.~3.1.8]{Zimmer84}), we hence have that $\uH$ is defined over
$\QQ$. The action of the $\QQ$-group $\uG$ on the $\QQ$-variety
$\uphahi$ is defined over $\QQ$, and the Galois group $\Gal(\CC|\QQ)$
acts on $\uphahi$ and on $\uG$ commuting with this action. For every
$\sigma\in \Gal(\CC|\QQ)$, we have $\uH^\sigma= \uH$. Hence by the
uniqueness of the point in $\phahi$ fixed by a conjugate of
$\PO(1,2)$, we have that $[Y]^\sigma=[Y]$ for every $\sigma\in
\Gal(\CC|\QQ)$. Thus $[Y]$ is a rational point. 
\end{proof}

An $\RR$-circle $C$ is {\em $K$-arithmetic} if its stabiliser in
$\Ga_K$ has a dense orbit in $C$. Proposition
\ref{prop:Fuchsianrationalpoint} explains this terminology: The
stabiliser in $\Ga_K$ of a $K$-arithmetic $\RR$-circle is arithmetic
(in the conjugate of $\PO(1,2)$ containing it). With $\omega=
\frac{-1+i\sqrt{3}}{2}$, the figure below shows part of the
$\Ga_{\QQ(\omega)}$-orbit of the standard infinite $\RR$-circle
$C_\infty$, which is $K$-arithmetic.

\medskip
\begin{center}
\includegraphics[width=9cm]{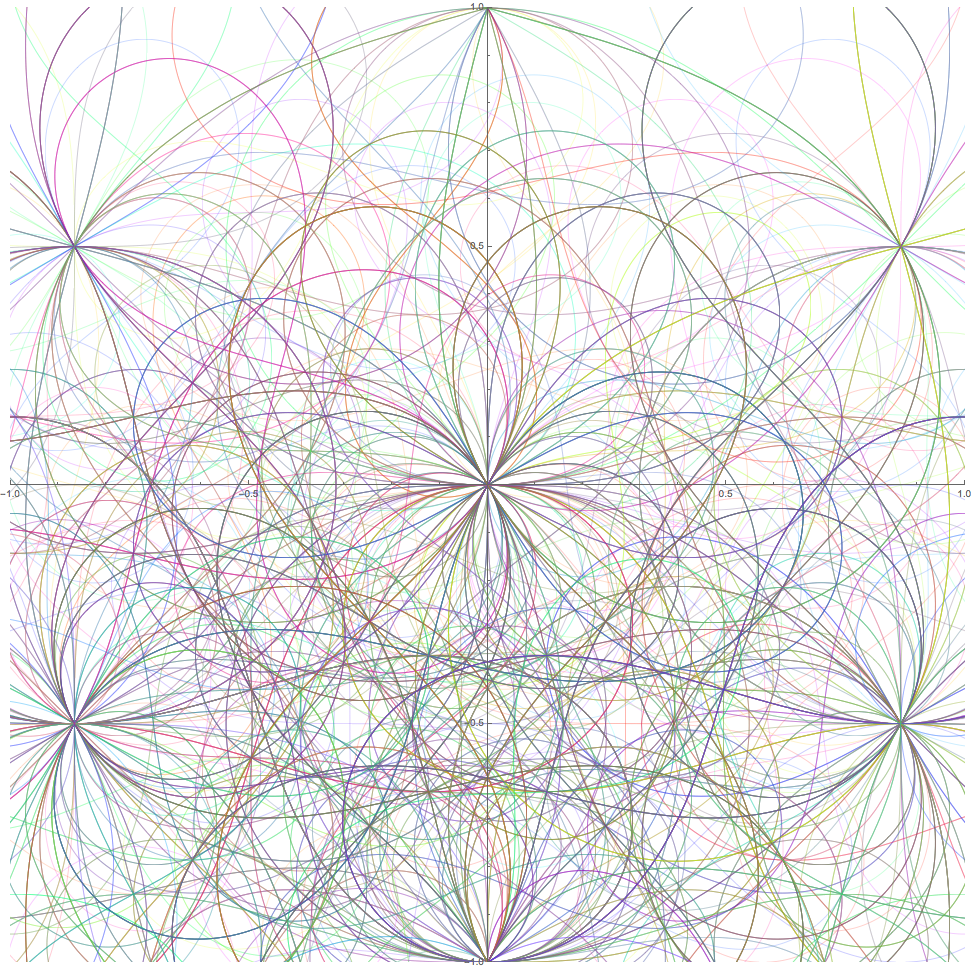}
\end{center}

\medskip
The next result reduces, up to commensurability and conjugacy 
in $\PU(1,2;K)$, the class of nonelementary $\RR$-Fuchsian subgroups that
we will study. Note that $\PU(1,2;K)$ is the commensurator of $\Ga_K$
in $\PU(1,2)$, see \cite[Theo.~2]{Borel66}.

\bprop \label{prop:reducpointrationel}
Any maximal nonelementary $\RR$-Fuchsian subgroup $\Ga$ of
$\Ga_K$ is commensurable up to conjugacy in $\PU(1,2;K)$ with the
stabiliser in $\Ga_K$ of the rational point $[Y_\Delta]\in\phahi$\footnote{or
equivalently to the stabiliser in $\Ga_K$ of the $\RR$-circle $C_\Delta$}
for some $\Delta\in\OOO_K$.  
If $\Delta\in \NN-\{0\}$ and if $$\ga_0=\begin{bmatrix} 
\frac{1+i}{2\sqrt{\Delta}} & 0 & \frac{1-i}{2\sqrt{\Delta}} \\ 
0 & 1 & 0 \\ 
\frac{(1-i)\sqrt{\Delta}}{2} & 0  & \frac{(1+i)\sqrt{\Delta}}{2}
\end{bmatrix}\,,$$ 
then $\ga_0\in\PU(1,2)$ and we have $\stab_{\Ga_K}[Y_\Delta]=
\ga_0\PO(1,2)\ga_0^{-1}\cap\Ga_K$.  
\eprop

\begin{proof} Let $\Ga$ be as in the statement. By Proposition
\ref{prop:Fuchsianrationalpoint}, there exists a rational point
$[Y]\in \uphahi(\QQ)$ in $\phahi$ such that $\Ga=\stab_{\Ga_K}[Y]=
\stab_{\Ga_K} C_Y$.  Up to conjugating $\Ga$ by an element in $\Ga_K$,
we may assume that the $\RR$-circle $C_Y$ is finite. The center of the
finite $\RR$-circle $C_Y$ belongs to $\PP_2(K)\cap (\partial_\infty
\HH^2_\CC-\{\infty\})$ by Lemma \ref{lem:calccentrerayon} (3). The
group of Heisenberg translations with coefficients in $K$ acts (simply
transitively) on $\PP_2(K)\cap(\partial_\infty\HH^2_\CC -\{\infty\})$.
Hence up to conjugating $\Ga$ by an element in $\PU(1,2;K)$, we may
assume that the center of the $\RR$-circle $C_Y$ is 
$0=[0:0:1]$.  By Lemma \ref{lem:formnormhahi} (and its proof), there
exists $\Delta\in K-\{0\}$ such that $[Y]=[Y_\Delta]$.  Since for
every $\lambda\in \CC^\times$ we have $\dil\lambda [Y_\Delta]\,
\ov{\dil\lambda}^{\,-1}= [Y_{\Delta\,\ov{\lambda}^{\,-2}}] $, up to
conjugating $\Ga$ by a Heisenberg dilation with coefficients in $K$,
we may assume that $\Delta\in\OOO_K$.

\smallskip
Fixing square roots of $\Delta$ and $\ov\Delta$ such that
$\sqrt{\ov\Delta}=\ov{\sqrt{\Delta}}$ , let 
$$
\ga'_0=\begin{bmatrix}
\frac{1+i}{2\sqrt{\ov{\Delta}}} & 0 & \frac{1-i}{2\sqrt{\ov{\Delta}}}
\\ 0 & 1 & 0 \\ \frac{(1-i)\sqrt{\Delta}}{2} & 0 & \frac{(1+i)
  \sqrt{\Delta}}{2} \end{bmatrix}\,.
$$ 
One easily checks using Equation \eqref{eq:equationsU} that
$\ga_0'\in\PU(1,2)$. An easy computation
proves that $\ga_0'[I_3]\,\ov{\ga_0'}^{\,-1}= \ga_0'\;\ov{\ga_0'}^{\,-1}
=[Y_\Delta]$.  Since the stabiliser of $[I_3]$ for the action of
$\PU(1,2)$ on $\phahi$ is equal to $\PO(1,2)$, the fact that
$$
\stab_{\Ga_K}[Y_\Delta] = \ga'_0\PO(1,2){\ga'_0}^{-1}\cap\Ga_K
$$ 
follows from the equivariance properties of the homeomorphisms in
Equation \eqref{eq:homeospaceRcercle}. Furthermore, $\ga'_0$ is the only
element of $\PU(1,2)$ satisfying this formula, up to right
multiplication by an element of $\PO(1,2)$. The last claim of
Proposition \ref{prop:reducpointrationel} follows since $\ga_0=\ga'_0$
when $\Delta\in\NN- \{0\}$.  
\end{proof}

Here is a geometric interpretation of the invariant $\Delta$
introduced in Proposition \ref{prop:reducpointrationel}: Since
$\ray(C_{Y_\Delta})= \ray(C_\Delta)= \frac{1}{\overline{\Delta}}$ for
every $\Delta\in\CC^\times$, the above proof shows that if the
$\RR$-circle $C_\Ga$ preserved by $\Ga$ is finite, then we may take
$\Delta\in\OOO_K-\{0\}$ squarefree (uniquely defined modulo a square
unit, hence uniquely defined if $D_K\neq -4,-3$) such that
$$
\Delta\in\;\left(\overline{\ray(C_\Ga)}\right)^{-1}\;(K^\times)^2\;.
$$

\section{Quaternion algebras, ternary quadratic forms and 
$\RR$-Fuchsian subgroups}
\label{sect:ternquad}

In this section, we describe the arithmetic structure of the maximal
nonelementary $\RR$-Fuchsian subgroups of $\Ga_K$.  By Proposition
\ref{prop:reducpointrationel}, it suffices to say from which
quaternion algebra the $\RR$-Fuchsian subgroup
$$
\Ga_{K,\,\Delta}=\stab_{\Ga_K}[Y_\Delta]
$$ 
arises for any $\Delta\in\OOO_K-\{0\}$.

Let $D,D'\in\QQ^\times$. The quaternion algebra $\Q=
\Big(\frac{D,D'}{\QQ}\Big)$ is the $4$-dimensional central simple
algebra over $\QQ$ with standard generators $i,j,k$ satisfying the
relations $i^2=D$, $j^2=D'$ and $ij=-ji=k$. If $x= x_0+x_1i+x_2j+x_3k$
is an element of $\Q$, we denote its {\it conjugate} by
$$
\ov x=x_0-x_1i-x_2j-x_3k\,,
$$ 
its (reduced) {\em trace} by
$$
\tr x=x+\ov x = 2x_0\,,
$$ 
and its {\em (reduced) norm} by
$$
\n(x_0+x_1i+x_2j+x_3k)=x\,\ov x=x_0^2-Dx_1^2-D'x_2^2+DD'x_3^2\;.
$$ 
The group of elements in $\Q(\ZZ)=\ZZ+i\ZZ+j\ZZ+k\ZZ$ with norm $1$
is denoted by $\Q(\ZZ)^1$. We refer to \cite{Vigneras80} and
\cite{MacRei03} for generalities on quaternion algebras.

The quaternion algebra $\Q$ {\em splits over $\RR$} if the
$\RR$-algebra $\Q(\RR)=\Q\otimes_\QQ\RR$ is isomorphic to the
$\RR$-algebra $\M_2(\RR)$ of $2$-by-$2$ matrices with real entries.
We say that a nonelementary $\RR$-Fuchsian subgroup $\Ga$ of
$\PU(1,2)$ {\it arises from} the quaternion algebra $\Q=
\big(\frac{D,D'}{\QQ}\big)$ if $\Q$ splits over $\RR$ and if there
exists a Lie group epimorphism $\varphi$ from $\Q(\RR)^1$ to the
conjugate of $\PO(1,2)$ containing $\Ga$, with kernel the center
$Z(\Q(\RR)^1)$ of $\Q(\RR)^1$, such that $\Ga$ and
$\varphi(\Q(\ZZ)^1)$ are commensurable.

Let $\A_\QQ$ be the set of isomorphism classes of quaternion
algebras over $\QQ$. For every $A\in \A_\QQ$, we denote by
$$
A_0=\{x\in A\;:\; \tr x=0\}
$$
the linear subspace of $A$ of {\it pure quaternions}, generated by
$i,j,k$.  Let $\T_\QQ$ be the set of isometry classes of nondegenerate
ternary quadratic forms over $\QQ$ with {\it
  discriminant}\footnote{the determinant of the associated matrix} a
square. It is well known (see for instance \cite[\S
  2.3--2.4]{MacRei03} and \cite[\S I.3]{Vigneras80}) that the map
$\Phi$ from $\A_\QQ$ to $\T_\QQ$, which associates to $A\in \A_\QQ$
the {\it restricted norm form} $\n_{\mid A_0}$, is a bijection.  The
map $\Phi$ has the following properties, for every $A\in \A_\QQ$.

\medskip
(1) If $a,b\in \QQ^\times$ and $A$ is (the isomorphism class of)
  $\big(\frac{a,b}{\QQ}\big)$, then $\Phi(A)$ is (the equivalence
  class of) $-a \,x_1^2 -b \,x_2^2 + ab \, x_3^2$, whose discriminant
  is $(ab)^2$.

\medskip
(2) If $a,b,c\in \QQ^\times$ with $abc$ a square in $\QQ$ and
  if $q\in \T_\QQ$ is (the equivalence class of) $-a \,x_1^2 -b
  \,x_2^2 + c \, x_3^2$, then $\Phi^{-1}(q)$ is (the isomorphism class
  of) $\big(\frac{a,b}{\QQ}\big)$, since if $abc= \lambda^2$ with
  $\lambda\in\QQ$, then the change of variables $(x'_1,x'_2,x'_3)=
  (x_1,x_2,\frac{\lambda}{ab}x_3)$ over $\QQ$ turns $q$ to the
  equivalent form $-a \,x_1^2 -b \,x_2^2 + ab \, x_3^2$.

\medskip
(3) The quaternion algebra $A$ splits over $\RR$ if and only if
  $\Phi(A)$ is isotropic over $\RR$ (that is, if the real quadratic
  form $\Phi(A)$ is indefinite), see \cite[Coro I.3.2]{Vigneras80}.

\medskip
(4) The map $\Theta_A$ from $A(\RR)^\times$ to the special
  orthogonal group $\SO_{\Phi(A)}$ of $\Phi(A)$, sending the class of
  an element $a$ in $A(\RR)^\times$ to the linear map $a_0\mapsto a
  a_0 a^{-1}$ from $A_0$ to itself, is a Lie group epimorphism with
  kernel the center of $A(\RR)^\times$ (see
  \cite[Th.~2.4.1]{MacRei03}).  If $A(\ZZ)=\ZZ+\ZZ i+\ZZ j+\ZZ k$ is
  the usual order in $A$, then $\Theta_A$ sends $A(\ZZ)^{1}$ to a
  subgroup commensurable with $\SO_{\Phi(A)}(\ZZ)$.

\medskip
\begin{proof}[Proof of Theorem \ref{theo:main}]
The set $P_\Delta$ of fixed points of the linear Hermitian anti-holomorphic
involution $z\mapsto Y_\Delta \,\ov{z}$ from $\CC^3$ to $\CC^3$ is a
real vector space of dimension $3$, equal to
$$
P_\Delta=\{z\in\CC^3\;:\; z=Y_\Delta \ov{z}\;\}=
\big\{(z_0,z_1,z_2)\in\CC^3\;:\; 
z_1=\ov{z_1},\;\; z_2=\Delta \ov{z_0}\;\big\}\;.
$$
Let $V$ be the vector space over $\QQ$ such that $V(\RR)=\CC^3$ and
$V(\QQ)=K^3$. Since the coefficients of the equations defining
$P_\Delta$ are in $\QQ$, there exists a vector subspace $W=W_\Delta$
of $V$ over $\QQ$ such that $W(\RR)=P_\Delta$.  The restriction to $W$
of the Hermitian form $h$, which is defined over $\QQ$, is a ternary
quadratic form $q=q_\Delta$ defined over $\QQ$, that we now compute.

Since $K=\QQ+i\sqrt{|D_K|}\,\QQ$, we write
$$
\Delta=u+i\sqrt{|D_K|}\,v
$$ 
with $u,v\in\QQ$, and the variables $z_j=x_j+i\sqrt{|D_K|}\,y_j$
with $x_j, y_j \in \RR$ for $j\in\{0,1,2\}$. If $(z_0, z_1,z_2)\in
P_\Delta$, we have
\begin{align*}
h(z_0, z_1,z_2)&=-\Re(z_2\ov{z_0})+|z_1|^2=
-\Re(\Delta\ov{z_0}^2)+|z_1|^2\\ &=
-u \,x_0^2+u|D_K|\,y_0^2-2|D_K|v\,x_0y_0+x_1^2\;.
\end{align*}
The right hand side of this formula is a ternary quadratic form
$q=q_\Delta$ on $P_\Delta$, whose coefficients are indeed in $\QQ$. It
is nondegenerate and has nonzero discriminant $-w$, where
$$
w=v^2D_K^2+u^2|D_K|= N(\Delta)|D_K|\in\QQ-\{0\}\,.
$$

By equivariance of the homeomorphisms in Equation
\eqref{eq:homeospaceRcercle} and as $\stab_{\PU(1,2)} [Y_\Delta]$ is
equal to its normaliser, the map from $\stab_{\PU(1,2)} [Y_\Delta]$ to
the projective orthogonal group $\PO_q$ of the quadratic space
$(P_\Delta,q)$, induced by the restriction map from
$\stab_{\operatorname{U}(1,2)} P_\Delta$ to $\operatorname{O}(q)$,
sending $g$ to $g_{\mid P_\Delta}$, is a Lie group isomorphism. It
sends the lattice $\Ga_{K,\,\Delta}$ to a subgroup commensurable with
the lattice $\PO_q(\ZZ)$ in $\PO_q$.  If we find a nondegenerate
quadratic form $q'=q'_\Delta$ equivalent to $q$ over $\QQ$ up to a
rational scalar multiple, whose discriminant is a rational square, and
which is isotropic over $\RR$, then $\Ga_{K,\,\Delta}$ arises from the
quaternion algebra $\Phi^{-1}(q')$, by Properties (3) and (4) of the
bijection $\Phi$.

First assume that $u=0$. By an easy computation, we have
$$
q=-\Big(- x_1^2-\frac{|D_K|v}{2}\,(x_0-y_0)^2+
\frac{|D_K|v}{2} \,(x_0+y_0)^2\Big)\;.
$$ 
The quadratic form $q'=- X_1^2 -\frac{|D_K|v}{2}\,X_2^2+
\frac{|D_K|v}{2}\, X_3^2$ over $\QQ$ is equivalent to $q$ over $\QQ$
up to sign. Its discriminant is the rational square
$(\frac{|D_K|v}{2})^2$, and $q'$ represents $0$ over $\RR$.  By
Property (2) of the bijection $\Phi$, we have $\Phi^{-1}(q')=
\big(\frac{1,\,\frac{|D_K|v}{2}}{\QQ}\big)=
\big(\frac{1,\,1}{\QQ}\big)$. Therefore if $u=0$, then
$\Ga_{K,\,\Delta}$ arises from the trivial quaternion algebra
$\M_2(\QQ)$.

\smallskip
Now assume that $u\neq 0$. By an easy computation, we have
\begin{align*}
q= & -\frac{1}{u}\big(- u\,x_1^2-(v^2D_K^2+u^2|D_K|)y_0^2+
(u\,x_0+|D_K|v\,y_0)^2\big)\\ =& 
-\frac{1}{u^2w}\big(- u^2w\,x_1^2-uw^2\,y_0^2+
uw(u\,x_0+|D_K|v\,y_0)^2\big) \;.
\end{align*} 
The quadratic form $q'=- uw^2\,X_1^2 -wu^2\,X_2^2+uw\,X_3^2$ is
equivalent to $q$ over $\QQ$ up to a scalar multiple in $\QQ$. Its
discriminant is the rational square $(uw)^4$ and it represents $0$
over $\RR$.  By Property (2) of the bijection $\Phi$, we have
$\Phi^{-1}(q')= \big(\frac{uw^2,\,wu^2}{\QQ}\big) =
\big(\frac{u,\,w}{\QQ}\big)$.  Therefore if $u\neq0$, since
$u=\frac{1}{2}\Tr\Delta$ and $w=N(\Delta)|D_K|$, then
$\Ga_{K,\,\Delta}$ arises from the quaternion algebra
$\big(\frac{2\Tr\Delta,\,N(\Delta)|D_K|}{\QQ}\big)$.  This concludes
the proof of Theorem \ref{theo:main}.  
\end{proof}

\bcoro \label{coro:uniqueness}
Let $\Delta,\Delta'\in \OOO_K-\{0\}$ with nonzero traces. The
maximal nonelementary $\RR$-Fuchsian subgroups $\Ga_{K,\,\Delta}$ and
$\Ga_{K,\,\Delta'}$ are commensurable up to conjugacy in $\PU(1,2)$ if
and only if the quaternion algebras $\big(\frac{2 \Tr\Delta ,\,
  N(\Delta) \,|D_K|}{\QQ}\big)$ and $\big(\frac{2 \Tr\Delta' ,\,
  N(\Delta') \,|D_K|}{\QQ}\big)$ over $\QQ$ are isomorphic.  
\ecoro

\begin{proof} 
Since the action of $\PU(1,2)$ on the set of $\RR$-planes $\P_\RR$ is
transitive, this follows from the fact that two arithmetic Fuchsian
groups are commensurable up to conjugacy in $\PSL_2(\RR)$ if and only
if their associated quaternion algebras are isomorphic (see
\cite{Takeuchi75}).  \end{proof}

To complement Theorem \ref{theo:main}, we give a more explicit version
of its proof in the special case when $\Delta \in\NN -\{0\}$.

\bprop \label{prop:quaternionalgebraarising}
Let $\Delta\in\NN-\{0\}$. The maximal nonelementary $\RR$-Fuchsian
subgroup $\Ga_{K,\,\Delta}$ arises from the quaternion algebra
$\big(\frac{\Delta,\,|D_K|}{\QQ}\big)$.  
\eprop

\begin{proof} Let $\Delta\in\NN-\{0\}$. Let $D=\frac{|D_K|}{4}$ if $D_K\equiv 0
\mod 4$ and $D=|D_K|$ otherwise, so that $\OOO_K\cap\RR=\ZZ$ and
$\OOO_K\cap i\RR= i\sqrt{D}\,\ZZ $. Let $D'= D\Delta$. We have
$D,D'\in\NN-\{0\}$. Let $\Q=\big(\frac{D,-D'}{\QQ}\big)$.

The matrices 
$$
e_0=\begin{pmatrix} 0 & -1 \\ 0 & \ \ 0 \end{pmatrix}, \ 
e_1=\begin{pmatrix} 1 & \ \ 0 \\ 0 & -1 \end{pmatrix}, \ 
e_2=\begin{pmatrix} 0 & 0 \\ 1 & 0 \end{pmatrix}
$$ 
form 
a basis of the Lie algebra $\sss\lll_2(\RR)=
\bigg\{\begin{pmatrix} x_1 & -x_0 \\ x_2 & -x_1 \end{pmatrix} \;:\;
x_0,x_1,x_2\in\RR\bigg\}$ of $\PSL_2(\RR)$. Note that
$$
-\det(x_0e_0+x_1e_1+x_2e_2)= - x_0x_2+x_1^2
$$ 
is the quadratic form restriction of $h$ to $\RR^3\subset \CC^3$.  We
thus have a well known {\it exceptional isomorphism} between
$\PSL_2(\RR)$ and the identity component $\SO_0(1,2)$ of
$\operatorname{O}(1,2)$, which associates to $g \in \PSL_2(\RR)$ the
matrix in the basis $(e_0,e_1,e_2)$ of the linear automorphism
$\Ad(g):X\mapsto gXg^{-1}$, which belongs to $\GL(\sss\lll_2(\RR))$.
We denote by $\Theta: \PSL_2(\RR) \ra \PU(1,2)$ the group isomorphism
onto its image $\PO(1,2)$ obtained by composing this exceptional
isomorphism first with the inclusion of $\SO_0(1,2)$ in
$\operatorname{U}(1,2)$, then with the canonical projection in
$\PU(1,2)$. Explicitly, we have by an easy computation
$$
\Theta: \begin{bmatrix} a & b \\ c & d \end{bmatrix} \mapsto
\begin{bmatrix} 
a^2 & 2ab & b^2 \\ ac & ad+bc & bd \\ c^2 & 2cd  & d^2 
\end{bmatrix}\;.
$$

We have a map $\sigma_{D,-D'}: \Q\to\mat_2(\RR)$ defined by
$$
(x_0+x_1i+x_2j+x_3k)\mapsto\begin{pmatrix}
x_0+x_1\sqrt D & (x_2+x_3\sqrt D)\sqrt{D'}\\
-(x_2-x_3\sqrt D)\sqrt{D'} & x_0-x_1\sqrt D
\end{pmatrix}\;.
$$ 
As is well-known\footnote{see for instance \cite{Katok92}}, the
induced map $\sigma: \Q(\RR)^1 \ra\PSLR$ is a Lie group
epimorphism with kernel $Z(\Q(\RR)^1)$, such that
$\sigma(\Q(\ZZ)^1)$ is a discrete subgroup of $\PSLR$.  With $\ga_0$
as in Proposition \ref{prop:reducpointrationel}, for all
$x_0,x_1,x_2,x_3\in\ZZ$, a computation gives that the element
$\ga_0\;\Theta\big(\sigma(x_0+x_1i+x_2j+x_3k)\big)\,\ga_0^{-1}$ of
$\PU(1,2)$ is equal to 
$$
\begin{bmatrix}
\vspace{.2cm}
a(x)&b(x)&c(x)/\Delta\\
\vspace{.2cm}
d(x)\sqrt\Delta&\n(x)& \overline{d(x)}\,/\sqrt\Delta\\
\overline{c(x)}\,\Delta&\overline{b(x)}&\overline{a(x)}
\end{bmatrix}\,,
$$
where
\begin{align*}
a(x) & =x_0^2 +Dx_1^2 + (2D'x_2x_3)i\sqrt{D}\,,\\
b(x) & =2(x_1x_2 + x_0x_3 + (x_1x_3 +\frac{x_0x_2}{D})i\sqrt{D})
\frac{\sqrt{DD'}}{\sqrt{\Delta}}\,,\\
c(x) & =DD'x_3^2 + D'x_2^2 + 2x_0x_1 \,i\sqrt{D}\,,\\
d(x) & =(x_0x_3 - x_1x_2 +(\frac{x_0x_2}{D}- x_1x_3)i\sqrt{D})\sqrt{DD'}\,.
\end{align*}

\bigskip
Let us consider the order $\OOO$ of $\Q$ defined by
$$
\OOO=\{x_0+x_1i+x_2j+x_3k\in \Q(\ZZ)\;:\; 
x_1,x_2,x_3\equiv 0 \mod D\}\;.
$$ 
Since $\frac{\sqrt{DD'}}{\sqrt{\Delta}}=D\in\ZZ$ and $\sqrt{DD'\Delta}
=D'\in\ZZ$, the above computation shows that the subgroup $\ga_0\;
\Theta(\sigma(\OOO^1))\,\ga_0^{-1}$ of $\PU(1,2)$ is contained in
$\Ga_K$. Since
$$
\Big(\frac{D,\,-D'}{\QQ}\Big)=
\Big(\frac{|D_K|,-|D_K|\Delta}{\QQ}\Big)=
\Big(\frac{|D_K|,\Delta}{\QQ}\Big)\;,
$$ 
the result follows.  
\end{proof}

\medskip\noindent{\bf Remark.} Note that by Hilbert's Theorem 90, if $\Delta'\in K$
satisfies $|\Delta'|=1$, then there exists $\Delta'' \in\OOO_K-\{0\}$
such that $\Delta'=\frac{\Delta''}{\ov{\Delta''}}$, so that the
Heisenberg dilation $\dil{{\Delta''}^{-1}}$ commensurates
$\Ga_{K,\,\Delta'}$ to $\Ga_{K,\,N(\Delta'')}$ and $N(\Delta'')$
belongs to $\NN-\{0\}$. Hence Proposition
\ref{prop:quaternionalgebraarising} implies that $\Ga_{K,\,\Delta'}$
arises from the quaternion algebra
$\big(\frac{N(\Delta''),\,|D_K|}{\QQ}\big)$.

\medskip

We conclude this paper by a series of arithmetic and geometric
consequences of the above determination of the quaternion algebras
associated with the maximal nonelementary $\RR$-Fuchsian subgroups of
the Picard modular groups. Their proofs follow closely the arguments
in \cite{Maclachlan86} pages 309 and 310, and a reader not interested
in the arithmetic details may simply admit that they follow by
formally replacing $-d$ by $d$ in the statements of loc.~cit.

Recall that given $a\in\ZZ-\{0\}$ and $p$ an odd positive prime not
dividing $a$, the {\it Legendre symbol} $\big(\frac{a}{b}\big)$ is
equal to $1$ if $a$ is a square mod $p$ and to $-1$
otherwise. Recall\footnote{See for instance \cite[page 91]{Samuel67}.}
that if $d\in\ZZ-\{0\}$ is squarefree, a positive prime $p$ is either
\begin{itemize}
\item[$\bullet$] {\it ramified} in $\QQ(\sqrt{d})$ when $p\divides d$
  if $p$ is odd, and when $d\equiv 2,3\;[4]$ if $p=2$,
\item[$\bullet$] {\it split} in $\QQ(\sqrt{d})$ when $p\notdivides d$
  and $\big(\frac{d}{p}\big)=1$ if $p$ is odd, and when $d\equiv
  1\;[8]$ if $p=2$,
\item[$\bullet$] {\it inert} in $\QQ(\sqrt{d})$ when $p\notdivides d$,
  and $\big(\frac{d}{p}\big)=-1$ if $p$ is odd, and when $d\equiv 5\;
  [8]$ if $p=2$.
\end{itemize}

\noindent 
Recall that a quaternion algebra $A$ over $\QQ$ is determined up to
isomorphism by the finite (with even cardinality) set
$\operatorname{RAM}(A)$ of the positive primes $p$ at which $A$ {\it
  ramifies}, that is, such that $A\otimes_\QQ\QQ_p$ is a division
algebra.

\bprop \label{prop:construcRfuchs} Let $A$ be an indefinite quaternion
algebra over $\QQ$. If the positive primes at which $A$ is ramified
are either ramified or inert in $\QQ(\sqrt{|D_K|})$, then there exists
a maximal nonelementary $\RR$-Fuchsian subgroup of $\Ga_K$ whose
associated quaternion algebra is $A$.  \eprop

\begin{proof}
Recall\footnote{See for instance \cite{Vigneras80}, in particular
  pages 32 and 37, and \cite[Chap.~III]{Serre70}.} that for all
$a,b\in\ZZ-\{0\}$ and for all positive primes $p$, the ($p$-){\it
  Hilbert symbol} $(a,b)_p$, equal to $-1$ if $\big(\frac{a,
  \,b}{\QQ_p}\big)$ is a division algebra and $1$ otherwise, is
symmetric in $a,b$, and satisfies $(a,bc)_p=(a,b)_p(a,c)_p$ and
\begin{equation}\label{eq:properthilbsymbmodp}
(a,b)_p=
\begin{cases} 
(-1)^{\frac{u-1}{2}\frac{v-1}{2}+\alpha\frac{v^2-1}{8}+\beta\frac{u^2-1}{8}}
& {\rm if}\;p= 2, a=2^\alpha u, b=2^\beta v, \;{\rm with}\; 
u,v \;{\rm odd}\\
\big(\frac{a}{p}\big) & 
{\rm if}\;p\neq 2, p\notdivides a, \;p\divides b, \;p^2\notdivides b\,.
\end{cases}
\end{equation}

Let $d=\frac{|D_K|}{4}$ if $D_K\equiv 0\;[4]$ and $d=|D_K|$ otherwise,
so that $d\in\NN-\{0\}$ is squarefree. Given $A$ as in the statement,
we may write $\operatorname{RAM}(A)=\{p_1,\dots,p_r,r_1,\dots,r_s\}$
with $p_i$ inert in $\QQ(\sqrt{d})$, and $r_i$ ramified in
$\QQ(\sqrt{d})$, so that the prime divisors of $d$ are $r_1,\cdots,
r_s, s_1,\cdots, s_k$, unless some $r_i$, say $r_1$, is equal to $2$
and $d\equiv 3\;[4]$, in which case the prime divisors of $d$ are
$r_2,\cdots, r_s,s_1,\cdots, s_k$.  As in \cite[page
  310]{Maclachlan86}, let $q$ be an odd prime different from all
$p_i,r_i,s_i$ such that

$\bullet$~ $q\equiv p_1\cdots p_r\;[8]$ if no $r_i$ is equal to $2$,
$q\equiv 5\,p_1\cdots p_r\;[8]$ if $r_i=2$ and $d\equiv 2\;[4]$ and 
$q\equiv 3\,p_1\cdots p_r\;[8]$ if $r_i=2$ and $d\equiv 3\;[4]$,

$\bullet$~ for every $i=1,\dots, s$, if $r_i$ is odd, then
$\big(\frac{q}{r_i}\big)= - \big(\frac{p_1\cdots\, p_r} {r_i}\big)$,

$\bullet$~ for every $i=1,\dots, k$, if $s_i$ is odd, then
$\big(\frac{q}{s_i}\big)= \big(\frac{p_1\cdots\, p_r} {s_i}\big)$.

With $\Delta=p_1\cdots p_r q$, which is a positive squarefree integer,
let us prove that $A$ is isomorphic to $\big(\frac{d, \,\Delta}
{\QQ}\big)$. This proves the result by Proposition
\ref{prop:quaternionalgebraarising}. By the characterisation of the
quaternion algebras over $\QQ$, we only have to prove that for every
positive prime $t$ not in $\operatorname{RAM}(A)$, we have
$(d,\Delta)_t=1$ and for every positive prime $t$ in
$\operatorname{RAM}(A)$, we have $(d,\Delta)_t=-1$. We distinguish in
the first case between $t=q$, $t=s_i$, $t=2$ and $t\neq q,s_1,\cdots,
s_k,2$, and in the second case between $t=p_i$ and $t=r_i$. By using
several times Equation \eqref{eq:properthilbsymbmodp} and the fact
that $\big(\frac{d}{q}\big)=1$ since $r+s$ is even as $A$ is
indefinite, the result follows (see the Appendix for details).
\end{proof}

 Recall that the {\it wide commensurability} class of a
subgroup $H$ of a given group $G$ is the set of subgroups of $G$
which are commensurable up to conjugacy to $H$. Two groups are {\it
  abstractly commensurable} if they have isomorphic finite index
subgroups.

\bcoro\label{coro:manyclasses} Every Picard modular group $\Ga_K$
contains infinitely many wide commensurability classes in $\PU(1,2)$
of (uniform) maximal nonelementary $\RR$-Fuchsian subgroups.  
\ecoro

Corollary \ref{coro:intro} of the introduction follows from Corollary
\ref{coro:manyclasses}. Note that there is only one wide commensurability
class of nonuniform maximal nonelementary
$\RR$-Fuchsian subgroups of $\Ga_K$, by \cite[Thm.~8.2.7]{MacRei03}.

\begin{proof} As seen in Corollary \ref{coro:uniqueness}, two maximal
nonelementary $\RR$-Fuchsian subgroups are commensurable up to
conjugacy in $\PU(1,2)$ if and only if their associated quaternion
algebras are isomorphic. Two such quaternion algebras are isomorphic
if and only if they ramify over the same set of primes. By Proposition
\ref{prop:construcRfuchs}, for every finite set $I$ with even cardinality of
positive primes which are inert over $\QQ(\sqrt{|D_K|})$, the quaternion
algebra with ramification set equal to $I$ is associated with a
maximal nonelementary $\RR$-Fuchsian subgroup. Since there are
infinitely many inert primes over $\QQ(\sqrt{|D_K|})$, the result follows.
\end{proof}

\bcoro 
Any arithmetic Fuchsian group whose associated quaternion algebra $A$ is
defined over $\QQ$ has a finite index subgroup isomorphic to an
$\RR$-Fuchsian subgroup of some Picard modular group $\Ga_K$.  
\ecoro

\begin{proof} As in \cite{Maclachlan86} page 310, if $\operatorname{RAM}(A)=
\{p_1,\cdots,p_n\}$, let $d\in\NN-\{0\}$ be such that
$\big(\frac{d}{p_i} \big) =-1$ if $p_i$ is odd and $d\equiv 5\;[8]$ if
$p_i=2$, so that $p_1,\dots,p_n$ are inert in $\QQ(\sqrt{d})$, and take 
$K=\QQ(\sqrt{-d})$.
\end{proof}

\bcoro For all quadratic imaginary number fields $K$ and $K'$, there
are infinitely many abstract commensurability classes of Fuchsian
subgroups with representatives in both Picard modular groups $\Ga_K$
and $\Ga_{K'}$.  
\ecoro

\begin{proof} There are infinitely many 
primes $p$ such that
$\big(\frac{|D_K|}{p} \big) = \big(\frac{|D_{K'}|}{p} \big) = -1$, hence
infinitely many finite subsets of them with an even number of elements.
\end{proof}

\section*{Acknowledgements}
The second author
  would like to thank the Isaac Newton Institute for Mathematical
  Sciences, Cambridge, for support and hospitality in April 2017
  during the programme ``Non-positive curvature group actions and
  cohomology''. This work was supported by EPSRC grant
  \textnumero\ EP/K032208/1 and by the French-Finnish CNRS grant PICS
  \textnumero\ 6950. We thank a lot John Parker, in particular for
  simplifying Lemma \ref{lem:calccentrerayon}, and Yves Benoist, whose
  suggestion to use the relation between quaternion algebras and
  ternary quadratic form was critical for the conclusion of Section
  \ref{sect:ternquad}. We also thank Gaëtan Chenevier for his help
  with the final corollaries.

\appendix
\section{Details on the proof of Proposition \ref{prop:construcRfuchs}}

Let us prove in preamble that
\begin{equation}\label{eq:doverq}
\big(\frac{d}{q}\big)=1\;.
\end{equation}

By using the quadratic reciprocity law for the Jacobi symbol, and its
multiplicativity properties, by the second and third assumptions on
$q$, since the $p_\ell$'s are inert in $\QQ(\sqrt{d})$, since $r+s$ is
even and $q-p_1\dots p_r\equiv 0\;[4]$ by the first assumption on $q$,
we have if the $r_j$'s, $s_i$'s and $p_\ell$'s are odd
\begin{align*}
\big(\frac{d}{q}\big)&=(-1)^{\frac{q-1}2\frac{d-1}2}
\big(\frac{q}{d}\big)=(-1)^{\frac{q-1}2\frac{d-1}2}
\prod_{j}\big(\frac{q}{r_j}\big)\prod_{i}
\big(\frac{q}{s_i}\big)\nonumber\\ &=(-1)^{s+\frac{q-1}2\frac{d-1}2}
\prod_{j}\big(\frac{p_1\dots p_r}{r_j}\big)\prod_{i}
\big(\frac{p_1\dots p_r}{s_i}\big)\nonumber\\ &=(-1)^{s+\frac{q-1}2\frac{d-1}2}
\big(\frac{p_1\dots p_r}{d}\big)=
(-1)^{s+(\frac{q-1}2-\frac{p_1\dots p_r-1}2)\frac{d-1}2}
\big(\frac{d}{p_1\dots p_r}\big)\nonumber\\&=
(-1)^{s+(\frac{q-p_1\dots p_r}2)\frac{d-1}2}
\prod_{\ell=1}^r\big(\frac{d}{p_\ell}\big)=
(-1)^{r+s+(\frac{q-p_1\dots p_r}2)\frac{d-1}2}=1\;.
\end{align*}

Recall that the Jacobi symbol satisfies $\big(\frac{2}{n}\big)
=(-1)^{\frac{n^2-1}{8}}$ for every odd positive integer $n$.
Similarly, if some $s_i$ is not odd (which implies that the $r_j$'s
and $p_\ell$'s are odd), say $s_1=2$, then with $d'=d/2$ which is odd,
since $q^2- (p_1\dots p_r)^2\equiv 0\;[16]$ by the first assumption on
$q$, we have
\begin{align*}
\big(\frac{d}{q}\big)&=\big(\frac{2}{q}\big)\big(\frac{d'}{q}\big)
=(-1)^{\frac{q^2-1}{8}+\frac{q-1}2\frac{d'-1}2}
\big(\frac{q}{d'}\big)=(-1)^{\frac{q^2-1}{8}+\frac{q-1}2\frac{d'-1}2}
\prod_{j}\big(\frac{q}{r_j}\big)\prod_{i\neq 1}
\big(\frac{q}{s_i}\big)\\ &=(-1)^{s+\frac{q^2-1}{8}+\frac{q-1}2\frac{d'-1}2}
\prod_{j}\big(\frac{p_1\dots p_r}{r_j}\big)\prod_{i\neq 1}
\big(\frac{p_1\dots p_r}{s_i}\big)\\ &=
(-1)^{s+\frac{q^2-1}{8}+\frac{q-1}2\frac{d'-1}2}
\big(\frac{p_1\dots p_r}{d'}\big)\\ &=
(-1)^{s+\frac{q^2-1}{8}+(\frac{q-1}2-\frac{p_1\dots p_r-1}2)\frac{d'-1}2}
\big(\frac{d'}{p_1\dots p_r}\big)\\&=
(-1)^{s+\frac{q^2-1}{8}+(\frac{q-1}2-\frac{p_1\dots p_r-1}2)\frac{d'-1}2}
\big(\frac{2}{p_1\dots p_r}\big)\big(\frac{d}{p_1\dots p_r}\big)\\&=
(-1)^{s+\frac{q^2-1}{8}+(\frac{q-1}2-\frac{p_1\dots p_r-1}2)\frac{d'-1}2-
\frac{(p_1\dots p_r)^2-1}{8}} \big(\frac{d}{p_1\dots p_r}\big)\\&=
(-1)^{s+(\frac{q-p_1\dots p_r}2)\frac{d'-1}2+\frac{q^2-(p_1\dots p_r)^2}{8}}
\prod_{\ell=1}^r\big(\frac{d}{p_\ell}\big)\\&=
(-1)^{r+s+(\frac{q-p_1\dots p_r}2)\frac{d'-1}2+\frac{q^2-(p_1\dots p_r)^2}{8}}=1\;.
\end{align*}

Similarly, if some $p_\ell$ is not odd (which implies that the $r_j$'s
and $s_i$'s are odd), say $p_1=2$, then $d\equiv 5\; [8]$ since $p_1$
is inert in $\QQ(\sqrt{d})$, hence $\frac{d-1}2\equiv 0\;[2]$ and
$\frac{d^2-1}8\equiv 1\;[2]$. Therefore
\begin{align*}
\big(\frac{d}{q}\big)&=(-1)^{\frac{q-1}2\frac{d-1}2}
\big(\frac{q}{d}\big)=(-1)^{\frac{q-1}2\frac{d-1}2}
\prod_{j}\big(\frac{q}{r_j}\big)\prod_{i}
\big(\frac{q}{s_i}\big)\nonumber\\ &=(-1)^{s+\frac{q-1}2\frac{d-1}2}
\prod_{j}\big(\frac{p_1\dots p_r}{r_j}\big)\prod_{i}
\big(\frac{p_1\dots p_r}{s_i}\big)\nonumber\\ &=(-1)^{s+\frac{q-1}2\frac{d-1}2}
\big(\frac{p_1\dots p_r}{d}\big)=(-1)^{s+\frac{q-1}2\frac{d-1}2}
\big(\frac{2}{d}\big)\big(\frac{p_2\dots p_r}{d}\big)\\ &=
(-1)^{s+\frac{d^2-1}{8}+(\frac{q-1}2-\frac{p_2\dots p_r-1}2)\frac{d-1}2}
\big(\frac{d}{p_2\dots p_r}\big)\nonumber\\&=
(-1)^{s+\frac{d^2-1}{8}+(\frac{q-p_2\dots p_r}2)\frac{d-1}2}
\prod_{\ell=2}^r\big(\frac{d}{p_\ell}\big)=
(-1)^{r-1+s+\frac{d^2-1}{8}+(\frac{q-p_2\dots p_r}2)\frac{d-1}2}=1\;.
\end{align*}

Similarly, assume that some $r_i$ is not odd (which implies that the
$s_i$'s and $p_\ell$'s are odd), say $r_1=2$. Since $r_1$ is ramified
in $\QQ(\sqrt{d})$, we have either $d\equiv 2\;[4]$ or $d\equiv
3\;[4]$. Assume first that $d\equiv 2\;[4]$. Then with $d'=d/2$ which
is odd, since the first assumption $q\equiv 5\,p_1\dots p_r\;[8]$ on
$q$ implies that $q- p_1\dots p_r\equiv 0\;[4]$ and that $\frac{q^2-
  (p_1\dots p_r)^2}8 \equiv 1\;[2]$ as the $p_\ell$'s are then odd, we
have
\begin{align*}
\big(\frac{d}{q}\big)&=\big(\frac{2}{q}\big)\big(\frac{d'}{q}\big)
=(-1)^{\frac{q^2-1}{8}+\frac{q-1}2\frac{d'-1}2}
\big(\frac{q}{d'}\big)=(-1)^{\frac{q^2-1}{8}+\frac{q-1}2\frac{d'-1}2}
\prod_{j\neq 1}\big(\frac{q}{r_j}\big)\prod_{i}
\big(\frac{q}{s_i}\big)\\ &=(-1)^{s-1+\frac{q^2-1}{8}+\frac{q-1}2\frac{d'-1}2}
\prod_{j\neq 1}\big(\frac{p_1\dots p_r}{r_j}\big)\prod_{i}
\big(\frac{p_1\dots p_r}{s_i}\big)\\ &=
(-1)^{s-1+\frac{q^2-1}{8}+\frac{q-1}2\frac{d'-1}2}
\big(\frac{p_1\dots p_r}{d'}\big)\\ &=
(-1)^{s-1+\frac{q^2-1}{8}+(\frac{q-1}2-\frac{p_1\dots p_r-1}2)\frac{d'-1}2}
\big(\frac{d'}{p_1\dots p_r}\big)\\&=
(-1)^{s-1+\frac{q^2-1}{8}+(\frac{q-1}2-\frac{p_1\dots p_r-1}2)\frac{d'-1}2}
\big(\frac{2}{p_1\dots p_r}\big)\big(\frac{d}{p_1\dots p_r}\big)\\&=
(-1)^{s-1+\frac{q^2-1}{8}+(\frac{q-1}2-\frac{p_1\dots p_r-1}2)\frac{d'-1}2-
\frac{(p_1\dots p_r)^2-1}{8}} \big(\frac{d}{p_1\dots p_r}\big)\\&=
(-1)^{s-1+(\frac{q-p_1\dots p_r}2)\frac{d'-1}2+\frac{q^2-(p_1\dots p_r)^2}{8}}
\prod_{\ell=1}^r\big(\frac{d}{p_\ell}\big)\\&=
(-1)^{r+s-1+(\frac{q-p_1\dots p_r}2)\frac{d'-1}2+\frac{q^2-(p_1\dots p_r)^2}{8}}=1\;.
\end{align*}
Assume secondly that $d\equiv 3\;[4]$, so that as already said we have
$d=r_2\dots r_ss_1\dots s_k$ which is odd. Then $\frac{d-1}2$ is odd
and $\frac{q-p_1\dots p_r}2$ is odd by the first assumption on $q$,
hence as above
\begin{align*}
\big(\frac{d}{q}\big)&=(-1)^{\frac{q-1}2\frac{d-1}2}
\big(\frac{q}{d}\big)=(-1)^{\frac{q-1}2\frac{d-1}2}
\prod_{j\neq 1}\big(\frac{q}{r_j}\big)\prod_{i}
\big(\frac{q}{s_i}\big)\\ &=(-1)^{s-1+\frac{q-1}2\frac{d-1}2}
\prod_{j\neq 1}\big(\frac{p_1\dots p_r}{r_j}\big)\prod_{i}
\big(\frac{p_1\dots p_r}{s_i}\big)\\ &=
(-1)^{s-1+\frac{q-1}2\frac{d-1}2}
\big(\frac{p_1\dots p_r}{d}\big)\\ &=
(-1)^{s-1+(\frac{q-1}2-\frac{p_1\dots p_r-1}2)\frac{d-1}2}
\big(\frac{d}{p_1\dots p_r}\big)\\&=
(-1)^{s-1+(\frac{q-p_1\dots p_r}2)\frac{d-1}2}
\prod_{\ell=1}^r\big(\frac{d}{p_\ell}\big)\\&=
(-1)^{r+s-1+(\frac{q-p_1\dots p_r}2)\frac{d-1}2}=1\;.
\end{align*}
This proves Equation \eqref{eq:doverq}.

\medskip
Let $t$ be a positive prime. First assume that $t$ does not belong to
$\operatorname{RAM}(A)$. Then one of the following case occurs: $t=q$,
$t=s_i\neq 2$ for some $i$, $t=s_i= 2$ for some $i$, $t=2\neq s_i$ for
every $i$, or $t\neq q,s_1,\cdots, s_k,2$.

If $t=q$, then $t\neq 2$, $t\notdivides d$, $t\divides \Delta$,
$t^2\notdivides \Delta$, so that by the second claim of Equation
\eqref{eq:properthilbsymbmodp} and by Equation \eqref{eq:doverq}, we
have as wanted
$$
(d,\Delta)_t=\big(\frac{d}{q}\big)=1\;.
$$ 

If $t= s_i\neq 2$ for some $i=1,\dots, k$, by the second claim of
Equation \eqref{eq:properthilbsymbmodp} since $s_i\divides d$,
$s_i^2\notdivides d$, and by the third assumption on $q$, we have as
wanted
$$
(d,\Delta)_t=(d,p_1\dots p_r)_{s_i}(d,q)_{s_i}=
\big(\frac{p_1\cdots\, p_r} {s_i}\big)\big(\frac{q}{s_i}\big)=
\big(\frac{q}{s_i}\big)^2=1\;.
$$ 

If $t= s_i = 2$ for some $i=1,\dots, k$, then $d$, which is
squarefree, is equal to $2d'$ for some odd $d'$.  The $p_i$'s are odd,
and $q\equiv p_1\cdots p_r\;[8]$ by the first assumption on $q$, hence
$\Delta$ is odd, of the form $(2\Delta'+1)^2+8\Delta''=1+8\Delta'''$
for some $\Delta', \Delta'', \Delta'''\in\NN$. By the first claim of
Equation \eqref{eq:properthilbsymbmodp}, we have as wanted
$$
(d,\Delta)_t=(-1)^{\frac{d'-1}{2}\frac{\Delta-1}{2}+ (\Delta^2-1)/8}=
(-1)^{\frac{d'-1}{2}(4\Delta''')+ 2\Delta'''+8{\Delta'''}^2}=1\;.
$$ 

If $t=2\neq s_i$ for every $i=1,\dots, k$, then $\Delta$ and $d$ are
odd, and $\Delta\equiv (p_1\dots p_r)^2 \;[8]$ by the first assumption
on $q$ so that $\Delta\equiv 1\; [4]$ since the $p_i$'s are odd. Hence by
the first claim of Equation \eqref{eq:properthilbsymbmodp}, we have as
wanted
$$
(d,\Delta)_t=(-1)^{\frac{d-1}{2}\frac{\Delta-1}{2}}=1\;.
$$ 

\medskip
Now assume that $t$ belongs to $\operatorname{RAM}(A)$. In particular
$t$ is equal either to $p_i$ for some $i=1,\dots, r$ or to $r_i$ for
some $i=1,\dots, s$.

If $t=p_i\neq 2$, then $t$ does not divide $d$ since $p_i$ is inert in
$\QQ(\sqrt{d})$, and $t$ divides $\Delta$ which is squarefree. Hence
by the second claim of Equation \eqref{eq:properthilbsymbmodp}, since
$p_i$ is inert in $\QQ(\sqrt{d})$, we have as wanted
$$
(d,\Delta)_t= \big(\frac{d}{p_i}\big)=-1\;.
$$ 

If $t=p_i= 2$, then $d\equiv 5\;[8]$ since $p_i$ is inert in
$\QQ(\sqrt{d})$, hence $d$ is odd, $\frac{d-1}{2}$ is even and
$\frac{d^2-1}8$ is odd. Moreover $\Delta=2\Delta'$ with $\Delta'$ odd.
Therefore by the first claim of Equation
\eqref{eq:properthilbsymbmodp}, we have as wanted
$$
(d,\Delta)_t= (-1)^{\frac{d-1}{2}\frac{\Delta'-1}{2}+ (d^2-1)/8}=-1\;.
$$ 

If $t=r_i\neq 2$ for some $i=1,\dots, k$, then $t$ divides $d$ which
is squarefree and does not divide $\Delta$, hence by the second
assumption on $q$, we have as wanted
$$
(d,\Delta)_t=(d,p_1\dots p_r)_{r_i}(d,q)_{r_i}= 
\big(\frac{p_1\dots p_r}{r_i}\big)\big(\frac{q}{r_i}\big)=-1\;.
$$ 

Assume at last that $t=r_i=2$ for some $i=1,\dots, s$. Since $r_i$ is
ramified in $\QQ(\sqrt{d})$, we have either $d\equiv 2\;[4]$ or
$d\equiv 3\;[4]$. If $d\equiv 2\;[4]$, then $t$ divides $d$ which is
equal to $2d'$ with $d'$ odd. But $t$ does not divide $\Delta$ which
is odd and, by the first assumption on $q$, we have $\Delta=
5(2\Delta'+1)^2+8\Delta''=5+8\Delta'''$ for some $\Delta', \Delta'',
\Delta'''\in\NN$. Hence by the first claim of Equation
\eqref{eq:properthilbsymbmodp}, we have as wanted
$$
(d,\Delta)_t=(-1)^{\frac{d'-1}{2}\frac{\Delta-1}{2}+ (\Delta^2-1)/8}=
(-1)^{\frac{d'-1}{2}(2+4\Delta''') + 3+2\Delta'''+8{\Delta'''}^2}=-1\;.
$$ 
If $d\equiv 3\;[4]$, then $d$ and $\Delta$ are odd and, by the first
assumption on $q$, we have $\Delta= 3(2\Delta'+1)^2 +8\Delta''\equiv
3\;[4]$. Hence by the first claim of Equation
\eqref{eq:properthilbsymbmodp}, we have as wanted
$$
(d,\Delta)_t=(-1)^{\frac{d-1}{2}\frac{\Delta-1}{2}}=-1\;.
$$ 

{\small \bibliography{lemniscate.bbl} }

\bigskip
{\small
\noindent \begin{tabular}{l} 
Department of Mathematics and Statistics, P.O. Box 35\\ 
40014 University of Jyv\"askyl\"a, FINLAND.\\
{\it e-mail: jouni.t.parkkonen@jyu.fi}
\end{tabular}
\medskip

\noindent \begin{tabular}{l}
Laboratoire de mathématique d'Orsay,\\
UMR 8628 Univ. Paris-Sud et CNRS,\\
Universit\'e Paris-Saclay,\\
91405 ORSAY Cedex, FRANCE\\
{\it e-mail: frederic.paulin@math.u-psud.fr}
\end{tabular}
}

\end{document}